\documentclass[12pt]{amsart}

\input xy	
\xyoption{all}

\newcommand{\uppermu}{{\overline{\mu}}}
\newcommand{\lowermu}{{\underline{\mu}}}
\newcommand{\homog}{{\operatorname{homog}}}
\newcommand{\length}{{\operatorname{length}}}
\newcommand{\BOX}{{\operatorname{Box}}}
\newcommand{\ev}{{\operatorname{ev}}}
\newcommand{\Prob}{{\operatorname{Prob}}}
\newcommand{\Ubar}{{\overline{U}}}
\newcommand{\fbar}{{\bar{f}}}
\newcommand{\Ohat}{{\hat{\OO}}}
\newcommand{\bad}{{\operatorname{bad}}}
\newcommand{\sing}{{\operatorname{sing}}}
\newcommand{\smooth}{{\operatorname{smooth}}}
\newcommand{\Xtwo}{X^{\{2\}}}
\newcommand{\Sigmatwo}{\Sigma^{\{2\}}}

\newcommand{\N}{{\mathbf N}}

\newcommand{\Q}{{\mathbf Q}}

\newcommand{\Fbar}{{\overline{\F}}}
\newcommand{\Z}{{\mathbf Z}}

\newcommand{\Aff}{{\mathbf A}}

\newcommand{\calE}{{\mathcal E}}
\newcommand{\calP}{{\mathcal P}}
\newcommand{\calQ}{{\mathcal Q}}
\newcommand{\calQmedium}{{\mathcal Q}^{\operatorname{medium}}}
\newcommand{\calQhigh}{{\mathcal Q}^{\operatorname{high}}}
\newcommand{\calQlarge}{{\mathcal Q}^{\operatorname{large}}}
\newcommand{\calR}{{\mathcal R}}
\newcommand{\calS}{{\mathcal S}}
\newcommand{\T}{{\mathcal T}}
\newcommand{\PP}{{\mathbf P}}
\newcommand{\PPhat}{{\hat{\mathbf P}}}
\newcommand{\OO}{{\mathcal O}}
\newcommand{\II}{{\mathcal I}}
\newcommand{\pp}{{\mathfrak p}}
\newcommand{\mm}{{\mathfrak m}}
\newcommand{\F}{{\mathbf F}}
\newcommand{\coker}{\operatorname{coker}}

\newcommand{\Alb}{\operatorname{Alb}}

\newcommand{\Spec}{\operatorname{Spec}}
\newcommand{\Proj}{\operatorname{Proj}}

\newcommand{\nichts}{{\ensuremath{\left.\right.}}}
\newcommand{\isom}{\simeq}
\newcommand{\del}{\partial}
\newcommand{\tensor}{\otimes}

\newcommand{\red}{{\operatorname{red}}}
\newcommand{\directsum}{\oplus}

\newtheorem{theorem}{Theorem}[section]
\newtheorem{lemma}[theorem]{Lemma}
\newtheorem{cor}[theorem]{Corollary}
\newtheorem{prop}[theorem]{Proposition}
\newtheorem{conj}[theorem]{Conjecture}

\theoremstyle{definition}
\newtheorem{example}[theorem]{Example}
\newtheorem{question}[theorem]{Question}

\theoremstyle{remark}
\newtheorem{rem}{Remark$\!\!$}	
\newtheorem{rems}{Remarks$\!\!$}	

\begin{document}

\title{Bertini theorems over finite fields}
\subjclass{Primary 14J70; Secondary 11M38, 11M41, 14G40, 14N05}
\author{Bjorn Poonen}
\thanks{This research was supported by NSF grant DMS-9801104
          and a Packard Fellowship.  Part of the research was done
	while the author was enjoying the hospitality of the
	Universit\'e de Paris-Sud.}
\address{Department of Mathematics, University of California, Berkeley, CA 94720-3840, USA}
\email{poonen@math.berkeley.edu}
\date{March 1, 2002}

\begin{abstract}
Let $X$ be a smooth quasiprojective subscheme of $\PP^n$ of
dimension $m \ge 0$ over $\F_q$.
Then there exist homogeneous polynomials $f$ over $\F_q$
for which the intersection of $X$ and the hypersurface $f=0$ is smooth.
In fact, the set of such $f$ has a positive density, 
equal to $\zeta_X(m+1)^{-1}$,
where $\zeta_X(s)=Z_X(q^{-s})$ is the zeta function of $X$.
An analogue for regular quasiprojective schemes over $\Z$
is proved, assuming the $abc$ conjecture and another conjecture.
\end{abstract}

\maketitle

\section{Introduction}
\label{introduction}

The classical Bertini theorems say that if a subscheme $X \subseteq \PP^n$
has a certain property, then for a sufficiently general hyperplane
$H \subset \PP^n$, $H \cap X$ has the property too.
For instance, if $X$ is a quasiprojective subscheme of $\PP^n$
that is smooth of dimension $m \ge 0$ over a field $k$,
and $U$ denotes the set of points $u$ in the dual projective
space $\PPhat^n$ corresponding to hyperplanes $H \subset \PP^n_{\kappa(u)}$
such that $H \cap X$ is smooth of dimension $m-1$ over the residue
field $\kappa(u)$ of $u$, then $U$ contains a dense open subset of $\PPhat^n$.
If $k$ is infinite, then $U \cap \PPhat^n(k)$ is nonempty, and hence one
can find $H$ over $k$.
But if $k$ is finite, then it can happen that
the finitely many hyperplanes $H$ over $k$
all fail to give a smooth intersection $H \cap X$.
See Theorem~\ref{ironic}.

N.~M.~Katz~\cite{katz1999} asked whether 
the Bertini theorem over finite fields can be salvaged by
allowing hypersurfaces of unbounded degree in place of hyperplanes.
(In fact he asked for a little more;
see Section~\ref{applicationsection} for details.)
We answer the question affirmatively below.
O.~Gabber~\cite[Corollary~1.6]{gabber2001} 
has independently proved the existence of good
hypersurfaces of any sufficiently large degree divisible by 
the characteristic of $k$.\phantom{\cite{katz2001}}

Let $\F_q$ denote the finite field of $q=p^a$ elements.
Let $S=\F_q[x_0,\dots,x_n]$ denote the homogeneous coordinate ring of $\PP^n$,
let $S_d \subset S$
denote the $\F_q$-subspace of homogeneous polynomials of degree $d$,
and let $S_\homog=\bigcup_{d=0}^\infty S_d$.
For each $f \in S_d$, let $H_f$ denote the subscheme
$\Proj(S/(f)) \subseteq \PP^n$.
Typically (but not always), $H_f$ is a hypersurface of dimension $n-1$
defined by the equation $f=0$.
Define the {\em density} of a subset $\calP \subseteq S_\homog$
by
	$$\mu(\calP):= \lim_{d \rightarrow \infty} 
		\frac{\#(\calP \cap S_d)}{\# S_d},$$
if the limit exists.
For a scheme $X$ of finite type over $\F_q$, 
define the zeta function~\cite{weil1949}
	$$\zeta_X(s)=Z_X(q^{-s}) 
	:= \prod_{\operatorname{closed }P \in X} 
			\left(1-q^{-s \deg P} \right)^{-1}
	= \exp \left( \sum_{r=1}^\infty 
			\frac{\#X(\F_{q^r})}{r} q^{-rs} \right).$$

\begin{theorem}[Bertini over finite fields]
\label{zeta}
Let $X$ be a smooth quasiprojective subscheme of $\PP^n$ of
dimension $m \ge 0$ over $\F_q$.
Define
	$$\calP :=\{\, f \in S_\homog: 
		H_f \cap X \text{ is smooth of dimension } m-1 \,\}.$$
Then $\mu(\calP)=\zeta_X(m + 1)^{-1}$.
\end{theorem}

\begin{rems}
$\nichts$
\begin{enumerate}
\item The empty scheme is smooth of any dimension, including $-1$.
\item In this paper, $\cap$ denotes scheme-theoretic intersection
(when applied to schemes).
\item If $n \ge 2$, the density is unchanged if we insist also that $H_f$
be a geometrically integral hypersurface of dimension $n-1$.
This follows from the easy Proposition~\ref{geometricallyintegral}.
\item The case $n=1$, $X=\Aff^1$, is a well known polynomial analogue of 
the fact that the set of squarefree integers
has density $\zeta(2)^{-1}=6/\pi^2$.
See Section~\ref{arithmetic} for a conjectural common generalization.
\item The density is independent of 
the choice of embedding $X \hookrightarrow \PP^n$!
\item By~\cite{dwork1960}, $\zeta_X$ is a rational function of $q^{-s}$,
so $\zeta_X(m+1)^{-1} \in \Q$.
\end{enumerate}
\end{rems}

The overall plan of the proof
is to start with all homogeneous polynomials of degree $d$, 
and then for each closed point $P \in X$
to sieve out the polynomials $f$ for which $H_f \cap X$ is singular at $P$.
The condition that $P$ be singular on $H_f \cap X$
amounts to $m+1$ linear conditions on the Taylor coefficients
of a dehomogenization of $f$ at $P$, and these linear conditions
are over the residue field of $P$.
Therefore one expects that the probability that $H_f \cap X$
is nonsingular at $P$ will be $1-q^{-(m+1)\deg P}$.
Assuming that these conditions at different $P$ are independent,
the probability that $H_f \cap X$ is nonsingular everywhere
should be
	$$\prod_{\operatorname{closed } P \in X} 
			\left(1-q^{-(m+1)\deg P} \right) 
	= \zeta_X(m+1)^{-1}.$$
Unfortunately, the independence assumption and 
the individual singularity probability estimates
break down once $\deg P$ becomes large relative to $d$.
Therefore we must approximate our answer by truncating the product
after finitely many terms, say those corresponding to $P$ of degree $<r$.
The main difficulty of the proof, 
as with many sieve proofs,
is in bounding the error of the approximation,
i.e., in showing that when $d \gg r \gg 1$,
the number of polynomials of degree $d$
sieved out by conditions
at the infinitely many $P$ of degree $\ge r$ is negligible.

In fact we will prove Theorem~\ref{zeta}
as a special case of the following,
which is more versatile in applications.
The effect of $T$ below is to prescribe the first few terms
of the Taylor expansions of the dehomogenizations of $f$ 
at finitely many closed points.

\begin{theorem}[Bertini with Taylor conditions]
\label{taylor}
Let $X$ be a quasiprojective subscheme of $\PP^n$ over $\F_q$.
Let $Z$ be a finite subscheme of $\PP^n$,
and assume that $U:=X-(Z \cap X)$ is smooth of dimension $m \ge 0$.
Fix a subset $T \subseteq H^0(Z,\OO_Z)$.
Given $f \in S_d$, let $f|_Z$ denote the element of $H^0(Z,\OO_Z)$
that on each connected component $Z_i$
equals the restriction
of $x_j^{-d} f$ to $Z_i$, 
where $j=j(i)$ is the smallest $j \in \{0,1\dots,n\}$
such that the coordinate $x_j$ is nonzero at $Z_i$.
Define
	$$\calP:=\{\, f \in S_\homog :
	\text{$H_f \cap U$ is smooth of dimension $m-1$, and } 
		f|_Z \in T \,\}.$$
Then
	$$\mu(\calP)=\frac{\#T}{\# H^0(Z,\OO_Z)} \; 
			\zeta_U(m + 1)^{-1}.$$
\end{theorem}

Using a formalism analogous to that of
Lemma~20 of~\cite{poonen-stoll1999},
we can deduce the following even stronger version,
which allows us to impose {\em infinitely} many local conditions,
provided that the conditions at most points are no more stringent
than the condition that 
the hypersurface intersect a given finite set of varieties smoothly.

\begin{theorem}[Infinitely many local conditions]
\label{formalism}
For each closed point $P$ of $\PP^n$ over $\F_q$,
let $\mu_P$ denote normalized Haar measure on the completed
local ring $\Ohat_P$ as an additive compact group,
and let $U_P$ be a subset of $\Ohat_P$ 
whose boundary $\del U_P$ has measure zero.
Also for each $P$, fix a nonvanishing coordinate $x_j$,
and for $f \in S_d$ let $f|_P$ denote the image of $x_j^{-d} f$ in $\Ohat_P$.
Assume that there exist smooth quasiprojective subschemes $X_1, \dots, X_u$
of $\PP^n$ of dimensions $m_i=\dim X_i$ over $\F_q$
such that
for all but finitely many $P$,
$U_P$ contains $f|_P$ whenever $f \in S_\homog$
is such that $H_f \cap X_i$ is smooth of dimension $m_i-1$ at $P$ for all $i$.
Define
	$$\calP:=\{\, f \in S_\homog :
	\text{$f|_P \in U_P$ for all closed points $P \in \PP^n$} \,\}.$$
Then 
	$\mu(\calP)=\displaystyle\prod_{\operatorname{closed } P \in \PP^n} 
		\mu_P(U_P)$.
\end{theorem}

\begin{rem}
Implicit in Theorem~\ref{formalism} is the claim
that the product $\prod_P \mu_P(U_P)$ always converges,
and in particular that its value is zero if and only if $\mu_P(U_P)=0$ 
for some closed point $P$.
\end{rem}

The proofs of Theorems \ref{zeta}, \ref{taylor}, and~\ref{formalism}
are contained in Section~\ref{bertiniproof}.
But the reader at this point is encouraged to jump
to Section~\ref{applicationsection} for applications,
and to glance at Section~\ref{arithmetic},
which shows that the $abc$ conjecture and another conjecture
imply analogues of our main theorems
for regular quasiprojective schemes over $\Spec \Z$.
The $abc$ conjecture is needed to apply a 
multivariable generalization~\cite{poonensquarefree}
of A. Granville's result~\cite{granville1998}
about squarefree values of polynomials.
For some open questions, 
see Sections \ref{openquestion} and~\ref{regularversussmooth},
and also Conjecture~\ref{unwanted}.

The author hopes that the technique of Section~\ref{bertiniproof} 
will prove useful in
removing the condition ``assume that the ground field $k$ is infinite''
from other theorems in the literature.

\section{Bertini over finite fields: the closed point sieve}
\label{bertiniproof}

Sections \ref{lowsection}, \ref{mediumsection}, and~\ref{highsection}
are devoted to the proofs of 
Lemmas \ref{lowdegree}, \ref{mediumdegree}, and~\ref{highdegree},
which are the main results needed in Section~\ref{proofs} to prove 
Theorems \ref{zeta}, \ref{taylor}, and~\ref{formalism}.

\subsection{Singular points of low degree}
\label{lowsection}

Let $A=\F_q[x_1,\dots,x_n]$ be the ring of regular functions on
the subset $\Aff^n:=\{x_0 \not=0\} \subseteq \PP^n$,
and identify $S_d$ with the set of dehomogenizations
$A_{\le d} = \{\, f\in A: \deg f \le d\,\}$,
where $\deg f$ denotes total degree.

\begin{lemma}
\label{surjective}
\nichts
\begin{enumerate}
\item[(a)] If $Y$ is a finite subscheme of $\PP^n$,
then the map 
$\phi_d: S_d=H^0(\PP^n,\OO_{\PP^n}(d)) \rightarrow H^0(Y, \OO_Y(d))$
is surjective for $d \gg 1$.
\item[(b)] If moreover $Y \subseteq \Aff^n := \{x_0 \not=0\}$,
then $\phi_d$ is surjective for $d \ge \dim H^0(Y,\OO_Y)-1$.
\end{enumerate}
\end{lemma}

\begin{proof}
(a) Let $\II_Y$ denote the ideal sheaf of $Y \subseteq \PP^n$.
Then $\coker(\phi_d)$ is contained in $H^1(\PP^n,\II_Y(d))$,
which vanishes for $d \gg 1$ by Theorem~III.5.2b of~\cite{hartshorne1977}.

\noindent(b) Dehomogenize by setting $x_0=1$,
so that $\phi_d$ is identified with a map from $A_{\le d}$
to $B:=H^0(Y,\OO_Y)$.
Let $b=\dim B$.
For $i \ge -1$, let $B_i$ denote the image of $A_{\le i}$ in $B$.
Then $0 = B_{-1} \subseteq B_0 \subseteq B_1 \subseteq \dots$,
so $B_j=B_{j+1}$ for some $j \in [-1,b-1]$.
Then 
	$$B_{j+2} = B_{j+1} + \sum_{i=1}^n x_i B_{j+1} 
		= B_j + \sum_{i=1}^n x_i B_j 
		= B_{j+1}.$$
Similarly $B_j=B_{j+1}=B_{j+2}=\dots$,
and these eventually equal $B$ by~(a).
Hence $\phi_d$ is surjective for $d \ge j$, 
and in particular for $d \ge b-1$.
\end{proof}

If $U$ is a scheme of finite type over $\F_q$,
let $U_{<r}$ denote the set of closed points of $U$ of degree $<r$.
Similarly define $U_{>r}$.

\begin{lemma}[Singularities of low degree]
\label{lowdegree}
Let notation and hypotheses be as in Theorem~\ref{taylor},
and define
	$$\calP_r:=\{\, f \in S_\homog: H_f \cap U 
		\text{ is smooth of dimension $m-1$ at all $P \in U_{<r}$, 
			and $f|_Z \in T$}\,\}.$$
Then
	$$\mu(\calP_r)= \frac{\#T}{\# H^0(Z,\OO_Z)} 
				\prod_{P \in U_{<r}}
				\left(1 - q^{-(m + 1)\deg P} \right).$$
\end{lemma}

\begin{proof}
Let $U_{<r}=\{P_1,\dots,P_s\}$.
Let $\mm_i$ be the ideal sheaf of $P_i$ on $U$,
let $Y_i$ denote the closed subscheme of $U$ corresponding to
the ideal sheaf $\mm_i^2 \subseteq \OO_U$,
and let $Y= \bigcup Y_i$.
Then $H_f \cap U$ is singular at $P_i$ (more precisely, not smooth
of dimension $m-1$ at $P_i$) if and only if
the restriction of $f$ to a section of $\OO_{Y_i}(d)$ is zero.
Hence $\calP_r \cap S_d$ is the inverse image of
$$T \times \prod_{i=1}^s \left(H^0(Y_i,\OO_{Y_i}) - \{0\} \right)$$
under the $\F_q$-linear composition
\begin{equation*}
	\phi_d: S_d=H^0(\PP^n,\OO_{\PP^n}(d)) \rightarrow 
		H^0(Y \cup Z, \OO_{Y \cup Z}(d)) 
		\isom H^0(Z,\OO_Z) \times 
		\prod_{i=1}^s H^0(Y_i,\OO_{Y_i}),
\end{equation*}
where the last isomorphism is the (noncanonical) untwisting, 
component by component, by division by the $d$-th powers
of various coordinates, as in the definition of $f|_Z$.
Applying part~(a) of Lemma~\ref{surjective} to $Y \cup Z$ 
shows that $\phi_d$ is surjective for $d \gg 1$, so
	$$\mu(\calP_r)=\lim_{d \rightarrow \infty}
			\frac{\# \left[ T \times \prod_{i=1}^s 
		\left(H^0(Y_i,\OO_{Y_i}) - \{0\} \right) \right]}
	{\# \left[ H^0(Z,\OO_Z) \times \prod_{i=1}^s 
		H^0(Y_i,\OO_{Y_i}) \right]}
	= \frac{\#T}{\# H^0(Z,\OO_Z)}
		\prod_{i=1}^s \left(1 - q^{-(m+1)\deg P_i} \right),$$
since $H^0(Y_i,\OO_{Y_i})$ has a two-step filtration whose quotients
$\OO_{U,P_i}/\mm_{U,P_i}$ and $\mm_{U,P_i}/\mm_{U,P_i}^2$
are vector spaces of dimensions $1$ and $m$ respectively
over the residue field of $P_i$.
\end{proof}

\subsection{Singular points of medium degree}
\label{mediumsection}

\begin{lemma}
\label{singularfraction}
Let $U$ be a smooth quasiprojective subscheme of $\PP^n$ of
dimension $m \ge 0$ over $\F_q$.
If $P \in U$ is a closed point of degree $e$,
where $e \le d/(m+1)$,
then the fraction of $f \in S_d$ such that $H_f \cap U$
is not smooth of dimension $m-1$ at $P$
equals $q^{-(m+1)e}$.
\end{lemma}

\begin{proof}
Let $\mm$ be the ideal sheaf of $P$ on $U$,
and let $Y$ denote the closed subscheme of $U$ corresponding to $\mm^2$.
The $f \in S_d$ to be counted are those in the kernel of
$\phi_d: H^0(\PP^n,\OO(d)) \rightarrow H^0(Y,\OO_Y(d))$.
We have $\dim H^0(Y,\OO_Y(d)) = \dim H^0(Y,\OO_Y) = (m+1)e$,
so $\phi_d$ is surjective by part~(b) of Lemma~\ref{surjective},
and the $\F_q$-codimension of $\ker \phi_d$ equals $(m+1)e$.
\end{proof}

Define the upper and lower densities 
$\uppermu(\calP)$, $\lowermu(\calP)$
of a subset $\calP \subseteq S$ as $\mu(\calP)$
was defined, but using $\limsup$ and $\liminf$ in place of $\lim$.

\begin{lemma}[Singularities of medium degree]
\label{mediumdegree}
Let $U$ be a smooth quasiprojective subscheme of $\PP^n$ of
dimension $m \ge 0$ over $\F_q$.
Define
\begin{align*}
	\calQmedium_r := \bigcup_{d \ge 0} \{\, f \in S_d: 
	& \;\; \text{there exists $P \in U$ 
				with $r \le \deg P \le \frac{d}{m+1}$}\\
	& \;\; \text{such that $H_f \cap U$ is not smooth 
		of dimension $m-1$ at $P$} \,\}.
\end{align*}
Then $\lim_{r \rightarrow \infty} \uppermu(\calQmedium_r)=0$.
\end{lemma}

\begin{proof}
Using Lemma~\ref{singularfraction} 
and the crude bound $\#U(\F_{q^e}) \le c q^{em}$ 
for some $c>0$ depending only on $U$ \cite{lang-weil1954},
we obtain
\begin{align*}
	\frac{\# (\calQmedium_r \cap S_d)}{\# S_d}
	&\le \sum_{e=r}^{\lfloor d/(m+1) \rfloor} 
		\text{ (number of points of degree $e$ in $U$) } q^{-(m+1)e} \\
	&\le \sum_{e=r}^{\lfloor d/(m+1) \rfloor} \#U(\F_{q^e}) q^{-(m+1)e} \\
	&\le \sum_{e=r}^\infty c q^{em} q^{-(m+1)e}, \\
	&= \frac{c q^{-r}}{1-q^{-1}}.
\end{align*}
Hence $\uppermu(\calQmedium_r) \le c q^{-r}/(1-q^{-1})$,
which tends to zero as $r\rightarrow \infty$.
\end{proof}

\subsection{Singular points of high degree}
\label{highsection}

\begin{lemma}
\label{singularfraction2}
Let $P$ be a closed point of degree $e$ in $\Aff^n$ over $\F_q$.
Then the fraction of $f \in A_{\le d}$ that vanish at $P$
is at most $q^{-\min(d,e^{1/n})}$.
\end{lemma}

\begin{proof}
Let $\nu=\min(d,e^{1/n})$, and let $\ev_P: A_{\le d} \rightarrow \F_{q^e}$
denote the evaluation-at-$P$ map.
If $e_i$ is the degree of the projection of $P$ onto the $i^{\text{th}}$ 
coordinate, then $e \le e_1 e_2 \cdots e_n$,
so some $e_i$ exceeds $e^{1/n}$.
For that $i$, $\ev_P$
{\em injects} $\sum_{j=0}^{\nu-1} \F_q x_i^j$ into $\F_{q^e}$.
Hence $\dim_{\F_q} \ev_P(A_{\le d}) \ge \nu$,
and the codimension of $\ker(\ev_P)$ in $A_{\le d}$ is at least $\nu$.
\end{proof}

\begin{lemma}[Singularities of high degree]
\label{highdegree}
Let $U$ be a smooth quasiprojective subscheme of $\PP^n$ of
dimension $m \ge 0$ over $\F_q$.
Define
	$$\calQhigh :=\bigcup_{d \ge 0} \{\, f \in S_d: 
		\exists P \in U_{>d/(m+1)}
	\text{ such that $H_f \cap U$ is not smooth 
		of dimension $m-1$ at $P$} \,\}.$$
Then $\uppermu(\calQhigh)=0$.
\end{lemma}

\begin{proof}
If the lemma holds for $U$ and for $V$, it holds for $U \cup V$,
so we may assume $U \subseteq \Aff^n$ is affine.

Given a closed point $u \in U$,
choose a system of local parameters $t_1,\dots,t_n \in A$
at $u$ on $\Aff^n$ such that $t_{m+1}=t_{m+2}=\dots=t_n=0$
defines $U$ locally at $u$.
Then $dt_1,\dots,dt_n$ are a $\OO_{\Aff^n,u}$-basis for 
the stalk $\Omega^1_{\Aff^n/\F_q,u}$.
Let $\del_1,\dots,\del_n$ be the dual basis of the stalk $\T_{\Aff^n/\F_q,u}$
of the tangent sheaf.
Choose $s \in A$ with $s(u) \not=0$ to clear denominators
so that $D_i:=s \del_i$ gives a global derivation $A \rightarrow A$
for $i=1,\dots,n$.
Then there is a neighborhood $N_u$ of $u$ in $\Aff^n$
such that $N_u \cap \{t_{m+1}=t_{m+2}=\dots=t_n=0\} = N_u \cap U$,
$\Omega^1_{N_u/\F_q}= \directsum_{i=1}^n \OO_{N_u} dt_i$,
and $s \in \OO(N_u)^*$.
We may cover $U$ with finitely many $N_u$,
so by the first sentence of this proof,
we may reduce to the case where $U \subseteq N_u$ for a single $u$.
For $f \in A_{\le d}$, $H_f \cap U$ 
{\em fails} to be smooth of dimension $m-1$ at a point $P \in U$ 
if and only if $f(P)=(D_1f)(P)=\dots=(D_mf)(P)=0$.

Now for the trick.
Let $\tau=\max_i(\deg t_i)$, $\gamma=\lfloor (d-\tau)/p \rfloor$,
and $\eta = \lfloor d/p \rfloor$.
If $f_0 \in A_{\le d}$, $g_1 \in A_{\le \gamma}$, \dots, 
$g_m \in A_{\le \gamma}$, 
and $h \in A_{\le \eta}$
are selected uniformly and independently at random,
then the distribution of 
	$$f:=f_0 + g_1^p t_1 + \dots + g_m^p t_m + h^p$$
is uniform over $A_{\le d}$.
We will bound the probability that an $f$ constructed in this way
has a point $P \in U_{>d/(m+1)}$ 
where $f(P)=(D_1f)(P)=\dots=(D_mf)(P)=0$.
By writing $f$ in this way, we partially decouple the $D_if$ from each other:
$D_if=(D_if_0) + g_i^p s$ for $i=1,\dots,m$.
We will select $f_0,g_1,\dots,g_m,h$ one at a time.
For $0 \le i \le m$, define $W_i:=U \cap \{D_1f=\dots=D_if=0\}$.

\bigskip
\noindent{\em Claim 1:} For $0 \le i \le m-1$,
conditioned on a choice of $f_0,g_1,\dots,g_i$
for which $\dim(W_i) \le m-i$,
the probability that $\dim(W_{i+1}) \le m-i-1$
is $1-o(1)$ as $d \rightarrow \infty$.
(The function of $d$ represented by the $o(1)$ depends on $U$ and the $D_i$.)

\bigskip
\noindent{\em Proof of Claim 1:}
Let $V_1$, \dots, $V_\ell$ denote the $(m-i)$-dimensional
$\F_q$-irreducible components of $(W_i)_\red$.
By B\'ezout's theorem~\cite[p.~10]{fulton1984},
	$$\ell \le (\deg \Ubar)(\deg D_1f)\dots(\deg D_if) = O(d^i)$$
as $d \rightarrow \infty$, where $\Ubar$ is the projective closure of $U$.
Since $\dim V_k \ge 1$,
there exists a coordinate $x_j$ depending on $k$ such that 
the projection $x_j(V_k)$ has dimension~1.
We need to bound the set
	$$G^\bad_k:=\{\,g_{i+1} \in A_{\le \gamma}: D_{i+1}f=(D_{i+1}f_0) 
		+ g_{i+1}^p s \text{ vanishes identically on $V_k$}\,\}.$$
If $g,g' \in G^\bad_k$,
then by taking the difference and multiplying by $s^{-1}$,
we see that $g-g'$ vanishes on $V_k$.
Hence if $G^\bad_k$ is nonempty, it is a coset of the subspace
of functions in $A_{\le \gamma}$ vanishing on $V_k$.
The codimension of that subspace, or equivalently
the dimension of the image of $A_{\le \gamma}$ in
the regular functions on $V_k$, exceeds $\gamma+1$,
since a nonzero polynomial in $x_j$ alone does not vanish on $V_k$.
Thus the probability that $D_{i+1}f$ vanishes on some $V_k$
is at most $\ell q^{-\gamma-1} = O(d^i q^{-(d-\tau)/p}) = o(1)$
as $d \rightarrow \infty$.
This proves Claim~1.

\bigskip
\noindent{\em Claim 2:} Conditioned on a choice of $f_0,g_1,\dots,g_m$
for which $W_m$ is finite,
$\Prob(H_f \cap W_m \cap U_{>d/(m+1)} = \emptyset) = 1-o(1)$
as $d \rightarrow \infty$.

\bigskip
\noindent{\em Proof of Claim 2:}
The B\'ezout theorem argument in the proof of Claim~1
shows that $\#W_m = O(d^m)$.
For a given point $P \in W_m$,
the set $H^\bad$ of $h \in A_{\le \eta}$ for which 
$H_f$ passes through $P$ is 
either $\emptyset$ 
or a coset of $\ker(\ev_P:A_{\le \eta} \rightarrow \kappa(P))$,
where $\kappa(P)$ is the residue field of $P$.
If moreover $\deg P>d/(m+1)$,
then Lemma~\ref{singularfraction2} implies 
$\#H^\bad / \# A_{\le \eta} \le q^{-\nu}$
where $\nu=\min \left(\eta,\left(\frac{d}{m+1} \right)^{1/n} \right)$.
Hence 
	$$\Prob(H_f \cap W_m \cap U_{>d/(m+1)} \not= \emptyset) \le
		\# W_m q^{-\nu} = O(d^m q^{-\nu}) = o(1)$$
as $d \rightarrow \infty$,
since $\nu$ grows like a positive power of $d$.
This proves Claim~2.

\bigskip
\noindent{\em End of proof of Lemma~\ref{highdegree}:}
Choose $f \in S_d$ uniformly at random.
Claims~1 and ~2 show that
with probability $\prod_{i=0}^{m-1} (1-o(1)) \cdot (1-o(1)) = 1-o(1)$
as $d \rightarrow \infty$,
$\dim W_i = m-i$ for $i=0,1,\dots,m$
and $H_f \cap W_m \cap U_{>d/(m+1)} = \emptyset$.
But $H_f \cap W_m$ is the subvariety of $U$
cut out by the equations $f(P)=(D_1f)(P)=\dots=(D_mf)(P)=0$,
so $H_f \cap W_m \cap U_{>d/(m+1)}$ is exactly the set
of points of $H_f \cap U$ of degree $>d/(m+1)$
where $H_f \cap U$ is not smooth of dimension $m-1$.
\end{proof}

\subsection{Proofs of theorems over finite fields}
\label{proofs}

\begin{proof}[Proof of Theorem~\ref{taylor}]
As mentioned in the proof of Lemma~\ref{mediumdegree},
the number of closed points of degree $r$ in $U$ is $O(q^{rm})$;
this guarantees that the product defining $\zeta_U(s)^{-1}$ 
converges at $s=m+1$.
By Lemma~\ref{lowdegree},
	$$\lim_{r \rightarrow \infty} \mu(\calP_r) 
	= \frac{\#T}{\# H^0(Z,\OO_Z)} \; \zeta_U(m+1)^{-1}.$$
On the other hand, the definitions imply
$\calP \subseteq \calP_r \subseteq \calP \cup \calQmedium_r \cup \calQhigh$,
so
$\uppermu(\calP)$ and $\lowermu(\calP)$
each differ from $\mu(\calP_r)$ by at most
$\uppermu(\calQmedium_r) + \uppermu(\calQhigh)$.
Applying Lemmas~\ref{mediumdegree} and~\ref{highdegree}
and letting $r$ tend to $\infty$, we obtain
	$$\mu(\calP)=\lim_{r \rightarrow \infty} \mu(\calP_r) 
	= \frac{\#T}{\# H^0(Z,\OO_Z)} \; \zeta_U(m+1)^{-1}.$$
\end{proof}

\begin{proof}[Proof of Theorem~\ref{zeta}]
Take $Z=\emptyset$ and $T=\{0\}$ in Theorem~\ref{taylor}.
\end{proof}

\begin{proof}[Proof of Theorem~\ref{formalism}]
The existence of $X_1,\dots,X_u$
and Lemmas~\ref{mediumdegree} and~\ref{highdegree}
let us approximate $\calP$
by the set $\calP_r$ defined only by the conditions 
at closed points $P$ of degree less than $r$,
for large $r$.
For each $P \in \PP^n_{<r}$,
the hypothesis $\mu_P(\del U_P)=0$
lets us approximate $U_P$ by a union of cosets of an ideal $I_P$ of
finite index in $\Ohat_P$.
(The details are completely analogous to those in the proof
of Lemma~20 of~\cite{poonen-stoll1999}.)
Finally, Lemma~\ref{surjective}(a)
implies that for $d \gg 1$, the images of $f \in S_d$
in $\prod_{P \in \PP^n_{<r}} \Ohat_P/I_P$
are equidistributed.
\end{proof}

Finally let us show that the densities
in our theorems
do not change if in the definition of density
we consider only $f$ for which $H_f$ is geometrically integral,
at least for $n \ge 2$.

\begin{prop}
\label{geometricallyintegral}
Suppose $n \ge 2$.
Let $\calR$ be the set of $f \in S_\homog$ for which $H_f$ fails to be 
a geometrically integral hypersurface of dimension $n-1$.
Then $\mu(\calR)=0$.
\end{prop}

\begin{proof}
We have $\calR = \calR_1 \cup \calR_2$
where $\calR_1$ is the set of $f \in S_\homog$
that factor nontrivially over $\F_q$,
and $\calR_2$ is the set of $f \in S_\homog$
of the form $N_{\F_{q^e}/\F_q}(g)$
for some homogeneous polynomial $g \in \F_{q^e}[x_0,\dots,x_n]$
and $e \ge 2$.
(Note: if our base field were not $\F_q$, an irreducible polynomial
that is not absolutely irreducible would be a constant times a norm,
but the constant is unnecessary here, 
since $N_{\F_{q^e}/\F_q}:\F_{q^e} \rightarrow \F_q$ is surjective.)

We have
$$\frac{\#(\calR_1 \cap S_d)}{\#S_d} \le
\frac{1}{\#S_d} \sum_{i=1}^{\lfloor d/2 \rfloor} (\# S_i)(\# S_{d-i}) 
		= \sum_{i=1}^{\lfloor d/2 \rfloor} q^{-N_i},$$
where 
$$N_i=\binom{n+d}{n} - \binom{n+i}{n} - \binom{n+d-i}{n}.$$
For $1 \le i \le d/2 - 1$,
\begin{align*}
	N_{i+1}-N_i &= \left[ \binom{n+d-i}{n} - \binom{n+d-i-1}{n} \right]
			- \left[ \binom{n+i+1}{n} - \binom{n+i}{n} \right] \\
	&= \binom{n+d-i-1}{n-1} - \binom{n+i}{n-1} \\
	&> 0.
\end{align*}
Similarly, for $d \gg n$,
$$N_1 = \binom{n+d-1}{n-1} - \binom{n+1}{n} 
	\ge \binom{n+d-1}{1} - \binom{n+1}{1} = d-2.$$
Thus 
$$\frac{\#(\calR_1 \cap S_d)}{\#S_d} \le
	\sum_{i=1}^{\lfloor d/2 \rfloor} q^{-N_i} 
	\le \sum_{i=1}^{\lfloor d/2 \rfloor} q^{2-d} \le d q^{2-d},$$
which tends to zero as $d \rightarrow \infty$.

The number of $f \in S_d$ that are norms
of homogeneous polynomials of degree $d/e$ over $\F_{q^e}$
is at most ${\left(q^e \right)}^{\binom{d/e + n}{n}}$.
Therefore
$$\frac{\#(\calR_2 \cap S_d)}{\#S_d} \le \sum_{e|d, e>1} q^{-M_e}$$
where $M_e = \binom{d+n}{n} - e \binom{d/e+n}{n}$.
For $2 \le e \le d$,
\begin{align*}
	\frac{e \binom{d/e+n}{n}} {\binom{d+n}{n}}
	&= \frac{e \left( \frac{d}{e} + n \right)
			\left( \frac{d}{e} + n -1 \right)
			\cdots
			\left( \frac{d}{e} + 1 \right)}
		{ (d+n)(d+n-1)\cdots(d+1) } \\
	&\le \frac{ e \left( \frac{d}{e} + n \right)
			\left( \frac{d}{e} + n -1 \right) }
		{ (d+n)(d+n-1) } \\
	&\le \frac{ e \left( \frac{d}{e} + n \right)^2}{d^2} \\
	&= \frac{1}{e} + \frac{2n}{d} + \frac{en^2}{d^2} \\
	&\le \frac{1}{2} + \frac{2n^2}{d} + \frac{dn^2}{d^2} \\
	&\le 2/3,
\end{align*}
once $d \ge 18n^2$.
Hence in this case, $M_e \ge \frac{1}{3} \binom{d+n}{n} \ge d^2/3$,
so 
$$\frac{\#(\calR_2 \cap S_d)}{\#S_d} \le \sum_{e|d, e>1} q^{-M_e}
	\le d q^{-d^2/3},$$
which tends to zero as $d \rightarrow \infty$.
\end{proof}

Another proof of Proposition~\ref{geometricallyintegral}
is given in Section~\ref{singularitysection},
but that proof is valid only for $n \ge 3$.

\section{Applications}
\label{applicationsection}

\subsection{Counterexamples to Bertini}

Ironically, we can use our 
hypersurface Bertini theorem to construct counterexamples
to the original hyperplane Bertini theorem!
More generally, we can show that hypersurfaces of bounded degree
do not suffice to yield a smooth intersection.

\begin{theorem}[Anti-Bertini theorem]
\label{ironic}
Given a finite field $\F_q$ and integers $n \ge 2$, $d \ge 1$,
there exists a smooth projective geometrically integral
hypersurface $X$ in $\PP^n$ over $\F_q$
such that
for each $f \in S_1 \cup \dots \cup S_d$,
$H_f \cap X$ fails to be smooth of dimension $n-2$.
\end{theorem}

\begin{proof}
Let $H^{(1)}$, \dots, $H^{(\ell)}$ be a list of the $H_f$ arising
from $f \in S_1 \cup \dots \cup S_d$.
For $i=1,\dots,\ell$ in turn, choose a closed point $P_i \in H^{(i)}$
distinct from $P_j$ for $j<i$.
Using a $T$ as in Theorem~\ref{taylor},
we can express the condition that a hypersurface in $\PP^n$ be smooth 
of dimension $n-1$ at $P_i$
and have tangent space at $P_i$ equal to that of $H^{(i)}$
whenever the latter is smooth of dimension $n-1$ at $P_i$.
Theorem~\ref{taylor} (with Proposition~\ref{geometricallyintegral})
implies that there exists a smooth projective geometrically integral
hypersurface $X \subseteq \PP^n$ satisfying these conditions.
Then for each $i$,
$X \cap H^{(i)}$ fails to be smooth of dimension $n-2$ at $P_i$.
\end{proof}

\begin{rem}
Katz~\cite[p.~621]{katz1999} remarks that if $X$ is the hypersurface 
	$$\sum_{i=1}^{n+1} (X_i Y_i^q - X_i^q Y_i) = 0$$
in $\PP^{2n+1}$ over $\F_q$ with homogeneous coordinates 
$X_1,\dots,X_{n+1},Y_1,\dots,Y_{n+1}$,
then $H \cap X$ is singular for every hyperplane $H$ in $\PP^{2n+1}$
over $\F_q$.
\end{rem}

\subsection{Singularities of positive dimension}
\label{singularitysection}

Let $X$ be a smooth quasiprojective subscheme of $\PP^n$ of
dimension $m \ge 0$ over $\F_q$.
Given $f \in S_\homog$, let $(H_f \cap X)_\sing$ denote
the closed subset of points 
where $H_f \cap X$ is not smooth of dimension $m-1$.

Although Theorem~\ref{zeta} shows that for a nonempty
smooth quasiprojective subscheme $X \subseteq \PP^n$
of dimension $m \ge 0$,
there is a positive probability that $(H_f \cap X)_\sing \not=\emptyset$,
we now show that the probability that $\dim (H_f \cap X)_\sing \ge 1$ 
is zero.

\begin{theorem}
\label{singularitylocus}
Let $X$ be a smooth quasiprojective subscheme of $\PP^n$ of
dimension $m \ge 0$ over $\F_q$.
Define
	$$\calS :=\{\, f \in S_\homog: \dim(H_f \cap X)_\sing \ge 1 \,\}.$$
Then $\mu(\calS)=0$.
\end{theorem}

\begin{proof}
This is a corollary of Lemma~\ref{highdegree} with $U=X$,
since $\calS \subseteq \calQhigh$.
\end{proof}

\begin{rem}
If $f \in S_\homog$ 
is such that $H_f$ is not geometrically integral of dimension $n-1$,
then $\dim (H_f)_\sing \ge n-2$.
Hence Theorem~\ref{singularitylocus} with $X=\PP^n$ gives a new proof 
of Proposition~\ref{geometricallyintegral}, at least when $n \ge 3$.
\end{rem}

\subsection{Space filling curves}

We next use Theorem~\ref{taylor} to answer affirmatively
all the open questions in~\cite{katz1999}.
In their strongest forms, these are
\begin{quote}
{\em Question 10:}
Given a smooth projective geometrically connected variety $X$
of dimension $m \ge 2$ over $\F_q$, and a finite extension $E$ of $\F_q$,
is there always a closed subscheme $Y$ in $X$, $Y \not=X$,
such that $Y(E)=X(E)$ and such that $Y$ is smooth and geometrically
connected over $\F_q$?

\medskip
\noindent
{\em Question 13:}
Given a closed subscheme $X \subseteq \PP^n$ over $\F_q$
that is smooth and geometrically connected of dimension $m$,
and a point $P \in X(\F_q)$,
is it true for all $d \gg 1$
that there exists a hypersurface $H \subseteq \PP^n$
of degree $d$ such that $P$ lies on $H$
and $H \cap X$ is smooth of dimension $m-1$?
\end{quote}
Both of these questions are answered by the following:

\begin{theorem}
\label{answersforkatz}
Let $X$ be a smooth quasiprojective subscheme of $\PP^n$ 
of dimension $m \ge 1$ over $\F_q$,
and let $F \subset X$ be a finite set of closed points.
Then there exists a smooth projective geometrically integral
hypersurface $H \subset \PP^n$ such that
$H \cap X$ is smooth of dimension $m-1$ and contains $F$.
\end{theorem}

\begin{rems}
\nichts
\begin{enumerate}
\item
If $m \ge 2$ and if $X$ in Theorem~\ref{answersforkatz}
is geometrically connected and projective in addition to being smooth,
then $H \cap X$ will be 
geometrically connected and projective too.
This follows from Corollary~III.7.9 in~\cite{hartshorne1977}.
\item
Recall that if a variety is geometrically connected and smooth,
then it is geometrically integral.
\item
Question~10 and (partially) Question~13 were independently
answered by Gabber~\cite{gabber2001}.
\end{enumerate}
\end{rems}

\begin{proof}[Proof of Theorem~\ref{answersforkatz}]
Let $T_{P,X}$ denote the Zariski tangent space of a point $P$ on $X$.
At each $P \in F$ choose a codimension~1 subspace $V_P \subset T_{P,\PP^n}$
not equal to $T_{P,X}$.
We will apply Theorem~\ref{formalism} with the following local
conditions: for $P \in F$, $U_P$ is the condition that
the hypersurface $H_f$ passes through $P$ and $T_{P,H}=V_P$;
for $P \not\in F$, $U_P$ is the condition that $H_f$
and $H_f \cap X$ be 
smooth of dimensions $n-1$ and $m-1$, respectively, at $P$.
Theorem~\ref{formalism} (with Proposition~\ref{geometricallyintegral})
implies the existence of a smooth projective geometrically integral
hypersurface $H \subset \PP^n$ satisfying these conditions.
\end{proof}

\begin{cor}
\label{everydimension}
Let $X$ be a smooth, projective, 
geometrically integral variety of dimension $m \ge 1$ over $\F_q$,
let $F$ be a finite set of closed points of $X$,
and let $y$ be an integer with $1 \le y \le m$.
Then there exists a smooth, projective, 
geometrically integral subvariety $Y \subseteq X$ of dimension $y$
such that $F \subset Y$.
\end{cor}

\begin{proof}
Use Theorem~\ref{answersforkatz} with reverse induction on $y$.
\end{proof}

\begin{cor}[Space filling curves]
\label{spacefillingcurves}
Let $X$ be a smooth, projective, 
geometrically integral variety of dimension $m \ge 1$ over $\F_q$,
and let $E$ be a finite extension of $\F_q$.
Then there exists a smooth, projective, 
geometrically integral curve $Y \subseteq X$ 
such that $Y(E)=X(E)$.
\end{cor}

\begin{proof}
Take $y=1$ and $F=X(E)$ in Corollary~\ref{everydimension}.
\end{proof}

In a similar way, we prove the following:
\begin{cor}[Space avoiding varieties]
Let $X$ be a smooth, projective, 
geometrically integral variety of dimension $m \ge 1$ over $\F_q$,
and let $\ell$ and $y$ be integers with $\ell \ge 1$ and $0 \le y \le m$.
Then there exists a smooth, projective, 
geometrically integral subvariety $Y \subseteq X$ of dimension $y$
such that $Y$ has no points of degree less than $\ell$.
\end{cor}

\begin{proof}
If $y=0$, let $Y$ be a closed point of $X$ of large degree.
(Such points exist since $X(\Fbar_q)$ is infinite.)
If $y>0$, repeat the arguments used in the proof of 
Theorem~\ref{answersforkatz} and Corollary~\ref{everydimension},
but in each application of Theorem~\ref{taylor},
force the hypersurface to avoid the finitely many points
of $X$ of degree less than $\ell$.
\end{proof}

\subsection{Albanese varieties}

For a smooth, projective, geometrically integral variety $X$ over a field,
let $\Alb X$ denote its Albanese variety.
As pointed out in~\cite{katz1999},
a positive answer to Question~13 implies that every positive dimensional
abelian variety $A$ over $\F_q$ contains a smooth, projective,
geometrically integral curve $Y$ such that the natural map
$\Alb Y \rightarrow A$ is surjective.
We generalize this slightly in the next result,
which strengthens Theorem~11 of~\cite{katz1999} in the finite field case.

\begin{theorem}
\label{albanese}
Let $X$ be a smooth, projective, 
geometrically integral variety of dimension $m \ge 1$ over $\F_q$.
Then there exists a smooth, projective, 
geometrically integral curve $Y \subseteq X$ 
such that the natural map $\Alb Y \rightarrow \Alb X$ is surjective.
\end{theorem}

\begin{proof}
Choose a prime $\ell$ not equal to the characteristic.
Represent each $\ell$-torsion point in $(\Alb X)(\Fbar_q)$
by a zero-cycle of degree zero on $X$,
and let $F$ be the finite set of closed points appearing in these.
Use Corollary~\ref{everydimension} to construct a smooth, projective,
geometrically integral curve $Y$ passing through all points of $F$.
The image of $\Alb Y \rightarrow \Alb X$ is an abelian subvariety
of $\Alb X$ containing all the $\ell$-torsion points,
so the image equals $\Alb X$.
(The trick of using the $\ell$-torsion points 
is due to Gabber~\cite{katz1999}.)
\end{proof}

\begin{rems}
\nichts
\begin{enumerate}
\item
A slightly more general argument
proves Theorem~\ref{albanese} 
over an arbitrary field $k$~\cite[Proposition~2.4]{gabber2001}.
\item
It is also true that any abelian variety over a field $k$
can be embedded as an abelian subvariety of the Jacobian of a
smooth, projective, geometrically integral curve 
over $k$~\cite{gabber2001}.
\end{enumerate}
\end{rems}

\subsection{Plane curves}

The probability that a projective plane curve over $\F_q$ is nonsingular
equals 
	$$\zeta_{\PP^2}(3)^{-1} = (1-q^{-1})(1-q^{-2})(1-q^{-3}).$$
(We interpret this probability as the density given by Theorem~\ref{zeta} 
for $X=\PP^2$ in $\PP^2$.)
Theorem~\ref{formalism} with a simple local calculation shows 
that the probability that
a projective plane curve over $\F_q$ has at worst nodes as singularities
equals 
	$$\zeta_{\PP^2}(4)^{-1} = (1-q^{-2})(1-q^{-3})(1-q^{-4}).$$
For $\F_2$, these probabilities are $21/64$ and $315/512$.

\begin{rem}
Although Theorem~\ref{zeta} guarantees the existence of
a smooth plane curve of degree $d$ over $\F_q$ only when
$d$ is sufficiently large relative to $q$,
in fact such a curve exists for every $d \ge 1$ 
and every finite field $\F_q$.
Moreover, the corresponding statement for hypersurfaces of specified
dimension and degree is true~\cite[\S11.4.6]{katz-sarnak1999}.
In fact, for any field $k$
and integers $n \ge 1$, $d \ge 3$ with $(n,d)$ not equal to $(1,3)$
or $(2,4)$, 
there exists a smooth hypersurface $X$ over $k$
of degree $d$ in $\PP^{n+1}$ such that
$X$ has {\em no nontrivial automorphisms} 
over $\overline{k}$~\cite{poonen2000}.
This last statement is false for $(1,3)$;
whether or not it holds for $(2,4)$ is an open question.
\end{rem}

\section{An open question}
\label{openquestion}

In response to Theorem~\ref{zeta}, 
Matt Baker has asked the following:

\begin{question}
\label{bakerquestion}
Fix a smooth quasiprojective subscheme $X$ of dimension $m$ over $\F_q$.
Does there exist an integer $n_0 > 0$
such that for $n \ge n_0$, if $\iota: X \rightarrow \PP^n$ is an embedding
and $\iota(X)$ is not contained in any hyperplane in $\PP^n$,
then there exists a hyperplane $H \subseteq \PP^n$ over $\F_q$
such that $H \cap \iota(X)$ is smooth of dimension $m-1$?
\end{question}

Theorem~\ref{zeta} proves that the answer is yes,
if one allows only the embeddings $\iota$ obtained by composing
a fixed initial embedding $X \rightarrow \PP^n$
with $d$-uple embeddings $\PP^n \rightarrow \PP^N$.
Nevertheless, we conjecture that for each $X$ of positive dimension,
the answer to Question~\ref{bakerquestion} is no.

\section{An arithmetic analogue}
\label{arithmetic}

We formulate an analogue of Theorem~\ref{zeta} in which the 
smooth quasiprojective scheme $X$ over $\F_q$
is replaced by a regular quasiprojective scheme $X$ over $\Spec \Z$,
and we seek hyperplane sections that are regular.
The reason for using regularity instead of the stronger condition
of being smooth over $\Z$ is discussed in Section~\ref{regularversussmooth}.

Fix $n \in \N = \Z_{\ge 0}$.
Redefine $S$ as the homogeneous coordinate ring $\Z[x_0,\dots,x_n]$ 
of $\PP^n_\Z$,
let $S_d \subset S$
denote the $\Z$-submodule of homogeneous polynomials of degree $d$,
and let $S_\homog=\bigcup_{d=0}^\infty S_d$.
If $p$ is prime, let $S_{d,p}$ denote the set of 
homogeneous polynomials in $\F_p[x_0,\dots,x_n]$ of degree $d$.
For each $f \in S_d$, let $H_f$ denote the subscheme
$\Proj(S/(f)) \subseteq \PP^n_\Z$.
Similarly, for $f \in S_{d,p}$, let $H_f$ denote
$\Proj(\F_p[x_0,\dots,x_n]/(f)) \subseteq \PP^n_{\F_p}$.

If $\calP$ is a subset of $\Z^N$ for some $N \ge 1$,
define the {\em upper density} 
	$$\uppermu(\calP) := \max_{\sigma} 
		\limsup_{B_{\sigma(1)} \rightarrow\infty}
		\cdots
		\limsup_{B_{\sigma(N)} \rightarrow\infty}
		\frac{\#(\calP \cap \BOX)}{\# \BOX},$$
where $\sigma$ ranges over permutations of $\{1,2,\ldots,N\}$
and 
	$$\BOX = \{(x_1,\dots,x_N) \in \Z^N : 
		|x_i| \le B_i \text{ for all $i$}\}.$$
(In other words, we take the $\limsup$ 
only over growing boxes whose dimensions
can be ordered so that each is very large relative to the previous dimensions.)
Define {\em lower density} $\lowermu(\calP)$ similarly
using $\min$ and $\liminf$.
Define upper and lower densities $\uppermu_d$ and $\lowermu_d$
of subsets of a fixed $S_d$
by identifying $S_d$ with $\Z^N$ using a $\Z$-basis of monomials.
If $\calP \subseteq S_\homog$, define 
$\uppermu(\calP) = \limsup_{d \rightarrow \infty} \uppermu_d(\calP \cap S_d)$
and 
$\lowermu(\calP) = \liminf_{d \rightarrow \infty} \lowermu_d(\calP \cap S_d)$.
Finally, if $\calP$ is a subset of $S_\homog$,
define $\mu(\calP)$ as the common value of $\uppermu(\calP)$
and $\lowermu(\calP)$ if $\uppermu(\calP)=\lowermu(\calP)$.
The reason for choosing this definition 
is that it makes our proof work;
aesthetically, we would have preferred to prove a stronger statement
by defining density as the limit over arbitrary boxes in $S_d$ 
with $\min\{d,B_1,\dots,B_N\} \rightarrow \infty$; 
probably such a statement is also true but extremely difficult to prove.

For a scheme $X$ of finite type over $\Z$, 
define the zeta function~\cite[\S1.3]{serre1965}
	$$\zeta_X(s) := \prod_{\operatorname{closed }P \in X} 
			\left(1-\#\kappa(P)^{-s} \right)^{-1},$$
where $\kappa(P)$ denotes the (finite) residue field of $P$.
This generalizes the definition of Section~\ref{introduction},
since a scheme of finite type over $\F_q$ can be viewed as
a scheme of finite type over $\Z$.
On the other hand, $\zeta_{\Spec \Z}(s)$ is the Riemann zeta function.

The $abc$ conjecture, formulated by D.~Masser and J.~Oesterl\'e
in response to insights of R.~C.~Mason, L.~Szpiro, and G.~Frey,
is the statement that for any $\epsilon > 0$,
there exists a constant $C=C(\epsilon)>0$ such that if $a,b,c$
are coprime positive integers satisfying $a+b=c$,
then $c < C ({\displaystyle \prod_{p | abc} p} )^{1+\epsilon}$.

For convenience, we say that a scheme $X$ of finite type over $\Z$
is {\em regular of dimension $m$} if for every {\em closed} point $P$ of $X$,
the local ring $\OO_{X,P}$ is regular of dimension $m$.
For a scheme $X$ of finite type over $\Z$,
this is equivalent to 
the condition that $\OO_{X,P}$ be regular for {\em all} $P \in X$
and all irreducible components of $X$ have Krull dimension $m$.
If $X$ is smooth of relative dimension $m-1$ over $\Spec \Z$,
then $X$ is regular of dimension $m$,
but the converse need not hold.
The empty scheme is regular of every dimension.

\begin{theorem}[Bertini for arithmetic schemes]
\label{arithmeticbertini}
Assume the $abc$ conjecture and Conjecture~\ref{unwanted} below.
Let $X$ be a quasiprojective subscheme of $\PP^n_\Z$ 
that is regular of dimension $m \ge 0$.
Define
	$$\calP :=\{\, f \in S_\homog: 
		H_f \cap X \text{ is regular of dimension } m-1 \,\}.$$
Then $\mu(\calP)=\zeta_X(m + 1)^{-1}$.
\end{theorem}

\begin{rem}
The case $X=\PP^0_\Z = \Spec \Z$ in $\PP^0_\Z$ 
of Theorem~\ref{arithmeticbertini}
is the statement that the density of squarefree integers is $\zeta(2)^{-1}$,
where $\zeta$ is the Riemann zeta function.
The proof of Theorem~\ref{arithmeticbertini} in general
will involve questions about squarefree values of multivariable polynomials.
\end{rem}

Given a scheme $X$, let $X_\Q = X \times \Q$,
and let $X_p = X \times \F_p$ for each prime $p$.

\begin{conj}
\label{unwanted}
Let $X$ be an integral quasiprojective subscheme of $\PP^n_\Z$
that is smooth over $\Z$ of relative dimension $r$.
There exists $c>0$ such that if $d$ and $p$ are sufficiently large, then
	$$\frac{\#\{\, f \in S_{d,p} : \dim (H_f \cap X_p)_\sing \ge 1 \,\}}
	{\# S_{d,p}} < \frac{c}{p^2}.$$
\end{conj}

Heuristically one expects that Conjecture~\ref{unwanted}
is true even if $c/p^2$ is replaced by $c/p^k$ for any fixed $k \ge 2$.
On the other hand, for the application to Theorem~\ref{arithmeticbertini},
it would suffice to prove a weak form of Conjecture~\ref{unwanted} 
with the upper bound $c/p^2$
replaced by any $\epsilon_p>0$ such that $\sum_p \epsilon_p <\infty$.
We used $c/p^2$ only to simplify the statement.

If $d$ is sufficiently large relative to $p$,
then Theorem~\ref{singularitylocus} provides a suitable upper bound
on the ratio in Conjecture~\ref{unwanted}.
If $p$ is sufficiently large relative to $d$,
then one can derive a suitable upper bound from the Weil Conjectures.
The difficulty lies in the case where $d$ is {\em comparable} to $p$.

\subsection{Singular points with small residue field}

We begin the proof of Theorem~\ref{arithmeticbertini}
with analogues of results in Section~\ref{lowsection}.
If $M$ is a finite abelian group, let $\length_\Z M$
denote its length as a $\Z$-module.

\begin{lemma}
\label{Zsurjective}
\nichts
\begin{enumerate}
\item[(a)] If $Y$ is a zero-dimensional closed subscheme of $\PP^n_\Z$,
then the map 
$\phi_d: S_d=H^0(\PP^n_\Z,\OO(d)) \rightarrow H^0(Y,\OO_Y(d))$
is surjective for $d \gg 1$.
\item[(b)] If moreover $Y \subseteq \Aff^n_\Z := \{x_0 \not=0\}$,
then $\phi_d$ is surjective for $d \ge \length_\Z H^0(Y,\OO_Y) - 1$.
\end{enumerate}
\end{lemma}

\begin{proof}
Copy the proof of Lemma~\ref{surjective}.
\end{proof}

\begin{lemma}
\label{cosets}
If $\calP \subseteq \Z^N$ is a union of $c$ distinct cosets of
a subgroup $G \subseteq \Z^N$ of index $m$,
then $\mu(\calP)=c/m$.
\end{lemma}

\begin{proof}
Without loss of generality, we may replace $G$ with its subgroup $(m\Z)^N$
of finite index.
The result follows,
since any of the boxes in the definition of $\mu$ can be approximated by 
a box of dimensions that are multiples of $m$, with an error that becomes
negligible compared with the number of lattice points in the box
as the box dimensions tend to infinity.
\end{proof}

If $X$ is a scheme of finite type over $\Z$,
define $X_{<r}$ as the set of closed points $P$ with $\#\kappa(P) < r$.
(This conflicts with the corresponding definition 
before Lemma~\ref{lowdegree}; forget that one.)
Define $X_{\ge r}$ similarly.
We say that $X$ is regular of dimension $m$ 
at a closed point $P$ of $\PP^n_\Z$
if either $P \not\in X$ or $\OO_{X,P}$ 
is a regular local ring of dimension $m$.

\begin{lemma}[Small singularities]
\label{smallfield}
Let $X$ be a quasiprojective subscheme of $\PP^n_\Z$ 
that is regular of dimension $m \ge 0$.
Define
	$$\calP_r:=\{\, f \in S_\homog: H_f \cap X 
		\text{ is regular of dimension $m-1$ at all $P \in X_{<r}$}
			   \,\}.$$
Then
	$$\mu(\calP_r)= \prod_{P \in X_{<r}} 
			\left(1 - \#\kappa(P)^{-(m + 1)} \right).$$
\end{lemma}

\begin{proof}
Given Lemmas \ref{Zsurjective} and~\ref{cosets},
the proof is the same as that of Lemma~\ref{lowdegree} with $Z=\emptyset$.
\end{proof}

\subsection{Reductions}
\label{reductionsection}

Theorem~1 of~\cite{serre1965}
shows that
$\displaystyle \prod_{P \in X_{<r}} \left(1 - \#\kappa(P)^{-(m + 1)} \right)$
converges to $\zeta_X(m+1)^{-1}$
as $r \rightarrow \infty$.
Thus Theorem~\ref{arithmeticbertini} follows from Lemma~\ref{smallfield}
and the following, whose proof will occupy the rest of 
Section~\ref{arithmetic}.

\begin{lemma}[Large singularities]
\label{largefield}
Assume the $abc$ conjecture and Conjecture~\ref{unwanted}.
Let $X$ be a quasiprojective subscheme of $\PP^n_\Z$ 
that is regular of dimension $m \ge 0$.
Define
\begin{align*}
	\calQlarge_r 
	:=\{\, f \in S_\homog: 
	& \text{ there exists $P \in X_{\ge r}$ such that } \\
	& \text{$H_f \cap X$ is not regular of dimension $m-1$ at $P$}
			   \,\}.
\end{align*}
Then $\lim_{r \rightarrow \infty} \uppermu(\calQlarge_r)= 0$.
\end{lemma}

Lemma~\ref{largefield} holds for $X$ if it holds for each subscheme
in an open cover of $X$, since by quasicompactness any such
open cover has a finite subcover.
In particular, we may assume that $X$ is connected.
Since $X$ is also regular, $X$ is integral.
If the image of $X \rightarrow \Spec \Z$ is a closed point $(p)$,
then $X$ is smooth of dimension $m$ over $\F_p$,
and Lemma~\ref{largefield} for $X$ follows from 
Lemmas \ref{mediumdegree} and~\ref{highdegree}.
Thus from now on, we assume that $X$ dominates $\Spec \Z$.

Since $X$ is regular, its generic fiber $X_\Q$ is regular.
Since $\Q$ is a perfect field, 
it follows that $X_\Q$ is smooth over $\Q$,
of dimension $m-1$.
By~\cite[17.7.11(iii)]{egaIV.III}, there exists an integer $t \ge 1$
such that $X \times \Z[1/t]$ is smooth of relative dimension $m-1$
over $\Z[1/t]$.

\subsection{Singular points of small residue characteristic}

\begin{lemma}[Singularities of small characteristic]
\label{smallcharacteristic}
Fix a nonzero prime $p \in \Z$.
Let $X$ be an integral quasiprojective subscheme of $\PP^n_\Z$ 
that dominates $\Spec \Z$ and is regular of dimension $m \ge 0$.
Define
\begin{align*}
	\calQ_{p,r} := \{\, f \in S_\homog: 
	& \;\;\text{there exists $P \in X_p$
		with $ \#\kappa(P) \ge r$}\\
	& \;\; \text{such that $H_f \cap X$ is not regular 
		of dimension $m-1$ at $P$}
			   \,\}.
\end{align*}
Then $\lim_{r \rightarrow \infty} \uppermu(\calQ_{p,r})= 0$.
\end{lemma}

\begin{proof}
We may assume that $X_p$ is nonempty.
Then, since $X_p$ is cut out in $X$ by a single equation $p=0$,
and since $p$ is neither a unit nor a zerodivisor in $H^0(X,\OO_X)$,
$\dim X_p=m-1$.

Let
\begin{align*}
	\calQmedium_{p,r} := \bigcup_{d \ge 0} \{\, f \in S_d: 
	& \;\; \text{there exists $P \in X_p$
			with $r \le \#\kappa(P) \le p^{d/(m+1)}$}\\
	& \;\; \text{such that $H_f \cap X$ is not regular 
		of dimension $m-1$ at $P$}
			   \,\}
\end{align*}
and
\begin{align*}
	\calQhigh_p := \bigcup_{d \ge 0} \{\, f \in S_d: 
	& \;\; \text{there exists $P \in X_p$
			with $ \#\kappa(P) > p^{d/(m+1)}$}\\
	& \;\; \text{such that $H_f \cap X$ is not regular 
		of dimension $m-1$ at $P$}
			   \,\}.
\end{align*}
Since $\calQ_{p,r} = \calQmedium_{p,r} \cup \calQhigh_p$, 
it suffices to prove
$\lim_{r \rightarrow \infty} \uppermu(\calQmedium_{p,r})=0$ 
and $\uppermu(\calQhigh_p)=0$.
We will adapt the proofs of
Lemmas \ref{mediumdegree} and~\ref{highdegree}.

If $P$ is a closed point of $X$,
let $\mm_{X,P} \subseteq \OO_X$ denote the ideal sheaf corresponding to $P$,
and let $Y_P$ be the closed subscheme of $X$ corresponding to 
the ideal sheaf $\mm_{X,P}^2$.
For fixed $d$, 
the set $\calQmedium_{p,r} \cap S_d$ is contained in the union over $P$
with $r \le \#\kappa(P) \le p^{d/(m+1)}$ 
of the kernel
of the restriction $\phi_P: S_d \rightarrow H^0(Y_P, \OO(d))$.
Since $H^0(Y_P,\OO(d)) \isom H^0(Y_P,\OO_{Y_P})$ 
has length $(m+1) [\kappa(P):\F_p] \le d$ as a $\Z$-module,
$\phi_P$ is surjective by Lemma~\ref{Zsurjective}(b),
and Lemma~\ref{cosets} implies $\uppermu(\ker \phi_P) = \#\kappa(P)^{-(m+1)}$.
Thus
	$$\uppermu(\calQmedium_{p,r} \cap S_d) 
	\le \sum_P \uppermu(\ker \phi_P) = \sum_P \#\kappa(P)^{-(m+1)}.$$
where the sum is over $P \in X_p$ 
with $r \le \#\kappa(P) \le p^{d/(m+1)}$.
The crude form $\# X_p(\F_{p^e}) = O(p^{e(m-1)})$ 
of the bound in~\cite{lang-weil1954} implies that
	$$\lim_{r \rightarrow \infty} \uppermu(\calQmedium_{p,r})
	= \lim_{r \rightarrow \infty} 
		\lim_{d \rightarrow \infty} 
			\uppermu(\calQmedium_{p,r} \cap S_d) = 0.$$

Next we turn to $\calQhigh_p$.
Since we are free to pass to an open cover of $X$,
we may assume that 
$X$ is contained in the subset 
$\Aff^n_\Z :=\{x_0 \not=0\}$ of $\PP^n_\Z$.
Let $A=\Z[x_1,\dots,x_n]$ be the ring of regular functions on $\Aff^n_\Z$.
Identify $S_d$ with the set of dehomogenizations
$A_{\le d} = \{\, f\in A: \deg f \le d\,\}$,
where $\deg f$ denotes total degree.

Let $\Omega=\Omega_{X_p/\F_p}$ denote the sheaf of regular
differentials on the reduced scheme associated to $X_p$.
For $P \in X_p$, define the dimension of the fiber
	$$\phi(P) = \dim_{\kappa(P)} \;
		\Omega \underset{\OO_{X_p}}\tensor \kappa(P).$$
Let $\mm_{X_p,P}$ be the maximal ideal of the local ring $\OO_{X_P,P}$.
If $P$ is a closed point of $X_p$, the isomorphism 
	$$\Omega \underset{\OO_{X_p}}\tensor \kappa(P) \isom 
		\frac{\mm_{X_p,P}}{\mm_{X_p,P}^2}$$
of Proposition~II.8.7 of~\cite{hartshorne1977}
shows that $\phi(P) = \dim_{\kappa(P)} \mm_{X_p,P}/\mm_{X_p,P}^2$;
moreover
	$$p \OO_{X,P} \rightarrow \frac{\mm_{X,P}}{\mm_{X,P}^2} 
	\rightarrow \frac{\mm_{X_p,P}}{\mm_{X_p,P}^2} \rightarrow 0$$
is exact.
Since $X$ is regular of dimension $m$, the middle term is 
a $\kappa(P)$-vector space of dimension $m$.
But the module on the left is generated by one element.
Hence $\phi(P)$ equals $m-1$ or $m$ at each closed point $P$.

Let $Y= \{\, P \in X_p : \phi(P) \ge m \,\}$.
By Exercise~II.5.8(a) of~\cite{hartshorne1977}, $Y$ is a closed subset,
and we give $Y$ the structure of a reduced subscheme of $X_p$.
Let $U = X_p - Y$.
Thus for closed points $P \in X_p$,
	$$\phi(P) = \begin{cases}
			m-1, & \text{if $P \in U$} \\
			m, & \text{if $P \in Y$.}
		\end{cases}$$
If $U$ is nonempty, then $\dim U = \dim X_p = m-1$,
so $U$ is smooth of dimension $m-1$ over $\F_p$,
and $\Omega|_U$ is locally free.
At a closed point $P \in U$,
we can find $t_1,\dots,t_n \in A$
such that $dt_1,\dots,dt_{m-1}$ represent 
an $\OO_{X_p,P}$-basis for the stalk $\Omega_P$,
and $dt_m,\dots,dt_n$ represent a basis for the kernel
of $\Omega_{\Aff^n/\F_p} \tensor \OO_{X_p,P} \rightarrow \Omega_P$.
Let $\del_1,\dots,\del_n \in \T_{\Aff^n/\F_p,P}$ 
be the basis of derivations dual to $dt_1,\dots,dt_n$.
Choose $s \in A$ nonvanishing at $P$ such that $s \del_i$
extends to a global derivation $D_i: A \rightarrow A$
for $i=1,2,\dots,m-1$.
In some neighborhood $V$ of $P$ in $\Aff^n_{\F_p}$, 
$dt_1,\dots,dt_n$ form a basis of $\Omega_{V/\F_p}$,
and $dt_1,\dots,dt_{m-1}$ form a basis of $\Omega_{U \cap V/\F_p}$,
and $s \in \OO(V)^*$.
By compactness, we may pass to an open cover of $X$ to assume $U \subseteq V$.
If $H_f \cap X$ is not regular at a closed point $Q \in U$,
then the image of $f$ in $\mm_{U,Q}/\mm_{U,Q}^2$ must be zero,
and it follows that $D_1 f$, \dots, $D_{m-1} f$, $f$ all vanish at $Q$.
The set of $f \in S_d$ such that there exists such a point in $U$
can be bounded using the induction argument
in the proof of Lemma~\ref{highdegree}.

It remains to bound the $f \in S_d$
such that $H_f \cap X$ is not regular at some closed point $P \in Y$.
Since $Y$ is reduced, 
and since the fibers of the coherent sheaf $\Omega \tensor \OO_Y$ on $Y$
all have dimension $m$,
Exercise~II.5.8(c) of~\cite{hartshorne1977} implies that
the sheaf is locally free.
By the same argument as in the preceding paragraph,
we can pass to an open cover of $X$,
and find $t_1,\dots,t_n,s \in A$
such that $dt_1,\dots,dt_n$ are a basis of the restriction of
$\Omega_{\Aff^n/\F_p}$ to a neighborhood of $Y$ in $\Aff^n_{\F_p}$,
and $dt_1,\dots,dt_m$ are an $\OO_Y$-basis of $\Omega \tensor \OO_Y$,
and $s \in \OO(Y)^*$ is such that if $\del_1,\dots,\del_n$ 
is the dual basis to $dt_1,\dots,dt_n$,
then $s \del_i$ extends to a derivation $D_i: A \rightarrow A$
for $i=1,\dots,m-1$.
(We could also define $D_i$ for $i=m$, but we already have enough.)
We finish again by using the induction argument
in the proof of Lemma~\ref{highdegree}.
\end{proof}

\subsection{Singular points of midsized residue characteristic}
\label{mediumcharacteristicsection}

While examining points of larger residue characteristic,
we may delete the fibers above small primes of $\Z$.
Hence in this section and the next, our lemmas
will suppose that $X$ is smooth over $\Z$.

\begin{lemma}[Singularities of midsized characteristic]
\label{midsizedcharacteristic}
Assume Conjecture~\ref{unwanted}.
Let $X$ be an integral quasiprojective subscheme of $\Aff^n_\Z$
that dominates $\Spec \Z$ 
and is smooth over $\Z$ of relative dimension $m-1$.
For $d,L,M \ge 1$, define
\begin{align*}
	\calQ_{d,L < \cdot < M} 
	:=\{\, f \in S_d: 
	& \text{ there exists $p \ge M$ and $P \in X_p$ such that } \\
	& \text{$H_f \cap X$ is not regular of dimension $m-1$ at $P$}
			   \,\}.
\end{align*}
Given $\epsilon>0$,
if $d$ and $L$ are sufficiently large,
then $\uppermu(\calQ_{d,L <\cdot < M}) < \epsilon$.
\end{lemma}

\begin{proof}
If $P$ is a closed point of degree at most $d/(m+1)$ over $\F_p$ 
where $L<p<M$,
then the set of $f \in S_d$ such that $H_f \cap X$ is not regular
of dimension $m-1$ at $P$ has upper density $\#\kappa(P)^{-(m+1)}$,
as in the argument for $\calQmedium_{p,r}$ in Lemma~\ref{smallcharacteristic}.
The sum over $\#\kappa(P)^{-(m+1)}$ over all such $P$
is small if $L$ is sufficiently large: 
this follows from~\cite{lang-weil1954}, as usual.
By Conjecture~\ref{unwanted}, 
the upper density of the set of $f \in S_d$ such that there 
exists $p$ with $L<p<M$ such that $\dim (H_f \cap X_p)_\sing \ge 1$
is bounded by $\sum_{L<p<M} c/p^2$, which again is small if $L$ is 
sufficiently large.

Let $\calE_{d,p}$ be the set of $f \in S_d$
for which $(H_f \cap X_p)_\sing$ is finite
and $H_f \cap X$ fails to be regular of dimension $m-1$
at some closed point $P \in X_p$ of degree greater than $d/(m+1)$ over $\F_p$.
It remains to show that if $d$ and $L$ are sufficiently large,
$\sum_{L<p<M} \uppermu(\calE_{d,p})$ is small.
Write $f=f_0+pf_1$ where $f_0$ has coefficients in $\{0,1,\dots,p-1\}$.
Once $f_0$ is fixed, $(H_f \cap X_p)_\sing$ is determined,
and in the case where it is finite, we let $P_1,\dots,P_\ell$ 
be its closed points of degree greater than $d/(m+1)$ over $\F_p$.
Now $H_f \cap X$ is not regular of dimension $m-1$ at $P_i$
if and only if the image of $f$ in $\OO_{X,P_i}/\mm_{X,P_i}^2$ is zero;
for fixed $f_0$, this is a condition only on the image
of $f_1$ in $\OO_{X_p,P_i}/\mm_{X_p,P_i}$.
It follows from Lemma~\ref{singularfraction2} that
the fraction of $f_1$ for this holds is at most $p^{-\nu}$
where $\nu=(d/(m+1))^{1/n}$.
Thus $\uppermu(\calE_{d,p}) \le \ell p^{-\nu}$
As usual, we may assume we have reduced to the case 
where $(H_f \cap X_p)_\sing$ is cut out by $D_1 f,\dots, D_{m-1} f,f$
for some derivations $D_i$, and hence by B\'ezout's Theorem, 
$\ell = O(d^m) = O(p^{\nu-2})$
as $d \rightarrow \infty$,
so $\uppermu(\calE_{d,p}) = O(p^{-2})$.
Hence $\sum_{L<p<M} \uppermu(\calE_{d,p})$ is small 
whenever $d$ and $L$ are large.
\end{proof}

\subsection{Singular points of large residue characteristic}
\label{largecharacteristicsection}

As in the previous section, 
$X$ denotes an integral quasiprojective subscheme of $\Aff^n_\Z$
that dominates $\Spec \Z$
and is smooth over $\Z$ of relative dimension $m-1$.
For fixed $d$, let $\Aff^N = \Aff^N_\Z$ denote the affine space
whose points correspond to polynomials of total degree at most $d$
in $x_0,\dots,x_n$.
Let $\Sigma \subseteq X \times \Aff^N$ denote the closed subscheme
of points $(x,f)$ such that the variety $H_f \cap X$ 
over the residue field of $(x,f)$ is not smooth of dimension $m-2$ at $x$.
Concretely, $\Sigma$ is the subscheme of $X \times \Aff^N$
locally cut out by the equations
$D_1 f = \cdots D_{m-1} f = f = 0$
where $f$ is written with indeterminate coefficients $c_i$
(which are the coordinates on $\Aff^N$)
and the $D_i$ are defined locally on $X$ as in the penultimate paragraph
of the proof of Lemma~\ref{smallcharacteristic}.

Let $\Xtwo = (X \times X) - \Delta$, where $\Delta$ is the image
of the diagonal map $X \rightarrow X \times X$.
Let $\Sigmatwo$ be the inverse image of $\Xtwo$
under the projection $\Sigma \times_{\Aff^N} \Sigma \rightarrow X \times X$.
Thus $\Sigmatwo$ is the closed subscheme of $\Xtwo \times \Aff^N$
whose points correspond to triples $(x_1,x_2,f)$
such that $x_1$ and $x_2$ are distinct points where $H_f \cap X$
fails to be smooth of dimension $m-2$.

In the following lemma, it is only the last part that will be used later.
\begin{lemma}
\label{sigmalemma}
If $d$ is sufficiently large, then:
\begin{enumerate}
\item[(a)] The projection $\Sigma \rightarrow X$ 
exhibits $\Sigma$ as a rank $N-m$ vector subbundle 
of $X \times \Aff^N \rightarrow X$.
\item[(b)] The projection $\Sigmatwo \rightarrow \Xtwo$ 
exhibits $\Sigmatwo$ as a rank $N-2m$ vector subbundle 
of $\Xtwo \times \Aff^N \rightarrow \Xtwo$.
\item[(c)] The image of the projection $\Sigmatwo_\Q \rightarrow \Aff^N_\Q$ 
is of dimension at most $N-2$ over $\Q$.
\item[(d)] The projection $\pi: \Sigma \rightarrow \Aff^N$ 
is a birational morphism onto its image $I$, 
and $I_\Q$ is an integral hypersurface in $\Aff^N_\Q$.
\item[(e)] There exists an integer $M>0$ and a squarefree polynomial 
$R(c_1,\dots,c_N) \in \Z[c_1,\ldots,c_N]$
such that if $\fbar$ is obtained from $f$ by specializing 
the coefficients $c_i$ to integers $\gamma_i$,
and if $H_\fbar \cap X$ fails to be regular 
at a closed point in the fiber $X_p$
for some prime $p \ge M$,
then $p^2$ divides the value $R(\gamma_1,\dots,\gamma_N)$.
\end{enumerate}
\end{lemma}

\begin{proof}
(a) The $m$ equations defining $\Sigma$ are linear
in the $c_i$ and are independent at each $x \in X$.

(b) Now there are $2m$ equations linear in the $c_i$.
If $d \ge 2m-1$ these are again independent,
by the argument of Lemma~\ref{surjective}(b).

(c) If $\Xtwo_\Q$ is nonempty, then by (b), 
	$$\dim \Sigmatwo_\Q = N-2m + \dim \Xtwo_\Q = N-2,$$
so the image of $\Sigmatwo_\Q$ in $\Aff^N_\Q$
has dimension at most $N-2$.

(d) By (a), $\Sigma$ is integral 
and smooth of relative dimension $N-1$ over $\Z$.
By definition of $I$, the variety $I_\Q$ is the image 
of $\Sigma_\Q \rightarrow \Aff^N_\Q$,
and the generic point of $\Sigma$ maps to the generic point $\eta$ of $I_\Q$.
Let $F$ be the fiber above of $\Sigma_\Q \rightarrow \Aff^N_\Q$ above $\eta$.
The fiber of $\Sigmatwo_\Q \rightarrow \Aff^N_\Q$
above $\eta$ equals $(F \times_{\eta} F) - \Delta_F$
where $\Delta_F$ is the diagonal in $F \times_{\eta} F$.
Hence if $F \rightarrow \eta$ is not an isomorphism
(that is, $\dim F > 0$, or $F$ has more than one {\em geometric} point),
then $\dim \Sigmatwo_\Q = 2 \dim F + \dim I_\Q$,
and we get the contradiction
	$$N-2 = \dim \Sigmatwo_\Q = 2 \dim F + \dim I_\Q
		\ge \dim F + \dim I_\Q = \dim \Sigma_\Q = N-1.$$
Thus $\Sigma_\Q \rightarrow I_\Q$ is birational.
In particular, $\dim I_\Q = \dim \Sigma_\Q = N-1$.
Moreover, since $\Sigma_\Q$ is integral, $I_\Q$ is integral as well.

(e)
Suppose that the equation defining $I_\Q$ in $\Aff^n_\Q$
is $R_0(c_1,\ldots,c_N)=0$.
By (d), $I_\Q$ is integral,
so we may assume that $R_0$ is an irreducible polynomial
in $\Z[c_1,\dots,c_N]$ with content $1$.
After inverting a finite number of nonzero primes of $\Z$,
we may assume that $R_0=0$ is also the equation defining
$I$ in $\Aff^n_\Z$.
Choose $M$ greater than all the inverted primes.
By (d), we may choose an open dense subset $I'$ of $I$ 
such that the birational morphism $\pi: \Sigma \rightarrow I$ induces
an isomorphism $\Sigma' \rightarrow I'$,
where $\Sigma'=\pi^{-1}(I')$.
By Hilbert's Nullstellensatz, there exists 
$R_1 \in \Z[c_1,\dots,c_N]$ such that $R_1$ vanishes on 
the closed subset $I-I'$ but not on $I$.
We may assume that $R_1$ is squarefree.
Define $R=R_0 R_1$.
Then $R$ is squarefree.

Suppose that $H_\fbar \cap X$ fails to be regular at a point $P \in X_p$
with $p \ge M$.
Let $\gamma$ denote the closed point of $\Aff^N$
defined by $c_1-\gamma_1=\cdots=c_N-\gamma_N=p=0$.
Then the point $(P,\gamma)$ of $X \times \Aff^N$ is in $\Sigma$.
Hence $\gamma \in I$, so $R_0(\gamma_1,\ldots,\gamma_N)$ is divisible by $p$.
If $\gamma \in I-I'$, 
then $R_1(\gamma_1,\ldots,\gamma_N)$ is divisible by $p$
as well, so $R(\gamma_1,\ldots,\gamma_N)$ is divisible by $p^2$,
as desired.

Therefore we assume from now on that $\gamma \in I'$, 
so $(P,\gamma) \in \Sigma'$.
Let $W$ be the inverse image of $I'$ under the closed immersion
$\Spec \Z \rightarrow \Aff^N$ defined by the ideal
$(c_1-\gamma_1,\ldots,c_N-\gamma_N)$.
Let $V$ be the inverse image of $\Sigma'$ under 
the morphism $X \hookrightarrow X \times \Aff^N$
induced by the previous closed immersion.
Thus we have a cube in which the top, bottom, front, and back faces are
cartesian:
$$\xymatrix{
V \ar[rr] \ar[dd] \ar[rd] & & {\Sigma'} \ar[dd] \ar[rd] \\
& X \ar[rr] \ar[dd] && X \times \Aff^N \ar[dd] \\
W \ar[rr] \ar[rd] && I' \ar[rd] \\
& {\Spec \Z} \ar[rr] && \Aff^N \\
}$$

Near $(P,\gamma) \in X \times \Aff^N$ the functions
$D_1 f, \cdots, D_{m-1} f, f $ cut out $\Sigma$
(and hence also its open subset $\Sigma'$)
locally in $X \times \Aff^N$.
Then $\OO_{V,P} = \OO_{X,P}/(D_1\fbar,D_2\fbar,\ldots,D_{m-1}\fbar,\fbar)$.
By assumption, $H_\fbar \cap X$ is not regular at $P$,
so $\fbar$ maps to zero in $\OO_{X,P}/\mm_{X,P}^2$.
Now $\OO_{X,P}$ is a regular local ring of dimension $m$,
$D_i \fbar \in \mm_{X,P}$, and $\fbar \in \mm_{X,P}^2$,
so the quotient $\OO_{V,P}$ has length at least $2$.
Since $\Sigma' \rightarrow I'$ is an isomorphism,
the cube shows that $V \rightarrow W$ is an isomorphism too.
Hence the localization of $W$ at $p$ has length at least $2$.
On the other hand $I'$ is an open subscheme of $I$,
whose ideal is generated by $R_0(c_1,\dots,c_N)$ 
(after some primes were inverted),
so $W$ is an open subscheme of $\Z/(R_0(\gamma_1,\ldots,\gamma_N))$.
Thus $R_0(\gamma_1,\ldots,\gamma_N)$ is divisible by $p^2$ at least.
Thus $R(\gamma_1,\ldots,\gamma_N)$ is divisible by $p^2$.
\end{proof}

Because of part~(e) of Lemma~\ref{sigmalemma},
we want to show that most values of a multivariable polynomial over $\Z$
are almost squarefree (that is, squarefree except for prime factors
less than $M$).
It is here that we need to use the $abc$ conjecture.

\begin{theorem}[Almost squarefree values of polynomials]
\label{squarefree}
Assume the $abc$ conjecture.
Let $F \in \Z[x_1,\dots,x_n]$ be squarefree.
For $M>0$, define 
	$$\calS_{M} := \{\, (a_1,\dots,a_n) \in \Z^n \mid
	F(a_1,\dots,a_n) 
	\text{ is divisible by $p^2$ for some prime $p \ge M$} \,\}.$$
Then $\uppermu(\calS_{M}) \rightarrow 0$ as $M \rightarrow \infty$.
\end{theorem}

\begin{proof}
The $n=1$ case is in~\cite{granville1998}.
The general case follows from Lemma~6.2 of~\cite{poonensquarefree},
in the same way that Corollary~3.3 there follows from Theorem~3.2 there.
Lemma~6.2 there is proved there by reduction to the $n=1$ case.
\end{proof}

\begin{rems}
\nichts
\begin{enumerate}
\item These results assume the $abc$ conjecture,
but the special case where $F$ factors into one-variable polynomials
of degree $\le 3$ is known unconditionally~\cite{hooley1967}.
Other unconditional results are contained in~\cite{gouvea-mazur1991}.
\item 
Theorem~\ref{squarefree} together with a simple sieve lets one
show that the naive heuristic
(multiplying probabilities for each prime $p$)
correctly predicts the density of $(a_1,\dots,a_n) \in \Z^n$
for which $F(a_1,\dots,a_n)$ is squarefree,
assuming the $abc$ conjecture.
\end{enumerate}
\end{rems}

\begin{lemma}[Singularities of large characteristic]
\label{largecharacteristic}
Assume the $abc$ conjecture.
Let $X$ be an integral quasiprojective subscheme of $\Aff^n_\Z$
that dominates $\Spec \Z$ 
and is smooth over $\Z$ of relative dimension $m-1$.
Define
\begin{align*}
	\calQ_{d,\ge M} 
	:=\{\, f \in S_d: 
	& \text{ there exists $p \ge M$ and $P \in X_p$ such that } \\
	& \text{$H_f \cap X$ is not regular of dimension $m-1$ at $P$}
			   \,\}.
\end{align*}
If $d$ is sufficiently large,
then $\lim_{M \rightarrow \infty} \uppermu(\calQ_{d,\ge M}) = 0$.
\end{lemma}

\begin{proof}
We may assume that $d$ is large enough for Lemma~\ref{sigmalemma}.
Apply Theorem~\ref{squarefree} to the squarefree polynomial $R$
provided by Lemma~\ref{sigmalemma}(e) for $X$.
\end{proof}

\subsection{End of proof}

We are now ready to prove Theorem~\ref{arithmeticbertini}.
Recall that in Section~\ref{reductionsection}
we reduced to the problem of proving Lemma~\ref{largefield}
in the case where $X$ is an integral quasiprojective 
subscheme of $\Aff^n_\Z$
such that $X$ dominates $\Spec \Z$ and is regular of dimension $m \ge 0$.
In Lemma~\ref{largefield}, $d$ tends to infinity for each fixed $r$,
and then $r$ tends to infinity.
We choose $L$ depending on $r$, and $M$ depending on $r$ and $d$,
such that $1 \ll L \ll r \ll d \ll M$.
(The precise requirement implied by each $\ll$ 
is whatever is needed below for the applications of the lemmas below.)
Then
\begin{equation}
\label{bigunion}
	\calQlarge_r \cap S_d \subseteq 
		\left( \bigcup_{p \le L} (\calQ_{p,r} \cap S_d)\right) 
		\cup \calQ_{d,L<\cdot<M} 
		\cup \calQ_{d,\ge M},
\end{equation}
and we will bound the upper density of each term on the right.
Recall from the end of Section~\ref{reductionsection} 
that $X$ has a subscheme of the form $X' = X \times \Spec \Z[1/t]$ 
that is smooth over $\Z$.
We may assume $L>t$.
By Lemma~\ref{smallcharacteristic}, 
$\lim_{r \rightarrow \infty} \uppermu(\calQ_{p,r}) = 0$ for each $p$,
so $\uppermu \left( \bigcup_{p\le L} \calQ_{p,r} \right)$ is small
(by which we mean tending to zero)
if $r$ sufficiently large relative to $L$.
By Lemma~\ref{midsizedcharacteristic} applied to $X'$,
if $L$ and $d$ are sufficiently large, then
$\uppermu(\calQ_{d,L<\cdot<M})$ is small.
By Lemma~\ref{largecharacteristic} applied to $X'$,
if $d$ is sufficiently large, and $M$ is sufficiently large relative to $d$,
then $\uppermu(\calQ_{d, \ge M})$ is small.
Thus by~(\ref{bigunion}), $\uppermu(\calQlarge_r)$ is small
whenever $r$ is large and $d$ is sufficiently large relative to $r$.
This completes the proof of Lemma~\ref{largefield} and hence
of Theorem~\ref{arithmeticbertini}. \hfill 

\begin{rem}
Arithmetic analogues of Theorems \ref{taylor} and~\ref{formalism},
and of many of the applications in Section~\ref{applicationsection}
can be proved as well.
\end{rem}

\subsection{Regular versus smooth}
\label{regularversussmooth}

One might ask what happens in Theorem~\ref{arithmeticbertini}
if we ask for $H_f \cap X$ to be not only regular,
but also smooth over $\Z$.
We now show unconditionally that this requirement is so strict,
that at most a density zero subset of polynomials $f$ satisfies it,
even if the original scheme $X$ is smooth over $\Z$.

\begin{theorem}
\label{smootharithmeticbertini}
Let $X$ be a nonempty quasiprojective subscheme of $\PP^n_\Z$ 
that is smooth of relative dimension $m \ge 0$ over $\Z$.
Define
	$$\calP^\smooth:=\{\, f \in S_\homog: 
		H_f \cap X \text{ is smooth 
			of relative dimension $m-1$ over $\Z$}\,\}.$$
Then $\mu(\calP^\smooth)=0$.
\end{theorem}

\begin{proof}
Let 
	$$\calP^\smooth_r :=\{\, f \in S_\homog: 
		H_f \cap X \text{ is smooth 
			of relative dimension $m-1$ over $\Z$
			at all $P \in X_{<r}$}\,\}.$$
Suppose $P \in X_{<r}$ lies above the prime $(p) \in \Spec \Z$.
Let $Y$ be the closed subscheme of $X_p$ corresponding to the ideal
sheaf $\mm^2$ where $\mm$ is the ideal sheaf of functions on $X_p$
vanishing at $P$.
Then for $f \in S_d$,
$H_f \cap X$ is smooth of relative dimension $m-1$ over $\Z$ at $P$
if and only if the image of $f$ in $H^0(Y,\OO(d))$ is nonzero.
Applying Lemma \ref{Zsurjective} to the union of such $Y$
over all $P \in X_{<r}$,
and using $\# H^0(Y,\OO(d)) = \#\kappa(P)^{m+1}$,
we find
	$$\mu(\calP^\smooth_r) = \prod_{P \in X_{<r}} 
			\left(1 - \#\kappa(P)^{-(m + 1)} \right).$$
Since $\dim X = m+1$, $\zeta_X(s)$ has a pole at $s=m+1$
and our product diverges to $0$ as $r \rightarrow \infty$.
(See Theorems 1 and 3(a) in~\cite{serre1965}.)
But $\calP^\smooth \subseteq \calP^\smooth_r$ for all $r$,
so $\mu(\calP^\smooth)=0$.
\end{proof}

A density zero subset of $S_\homog$ can still be nonempty
or even infinite.
For example, if $X=\Spec \Z[1/2,x] \hookrightarrow \PP^1_\Z$,
then $\calP^\smooth \cap S_d$ is infinite for infinitely many $d$:
$H_f \cap X$ is smooth over $\Z$ whenever $f$ is the homogenization
of $(x-a)^{2^b}-2$ for some $a,b \in \Z$ with $b \ge 0$.

On the other hand, N. Fakhruddin has given the following two examples
in which $\calP^\smooth \cap S_d$ is empty for all $d>0$.

\begin{example}
Let $X$ be the image of the $4$-uple embedding $\PP^1_\Z \rightarrow \PP^4_\Z$.
Then $X$ is smooth over $\Z$.
If $f \in \calP^\smooth \cap S_d$ for some $d>0$,
then $H_f \cap X \isom \coprod \Spec A_i$
where each $A_i$ is the ring of integers of a number field $K_i$
unramified above all finite primes of $\Z$,
such that $\sum [K_i:\Q] = 4d$.
The only absolutely unramified number field is $\Q$,
so each $A_i$ is $\Z$, and $H_f \cap X \isom \coprod_{i=1}^{4d} \Spec \Z$.
Then $4d = \# (H_f \cap X)(\F_2) \le \#X(\F_2) = \# \PP^1(\F_2) = 3$,
a contradiction.
\end{example}

\begin{example}
Let $X$ be the image of the $3$-uple embedding $\PP^2_\Z \rightarrow \PP^9_\Z$.
Then $X$ is smooth over $\Z$.
If $f \in \calP^\smooth \cap S_d$ for some $d>0$,
then $H_f \cap X$ is isomorphic to 
a smooth proper geometrically connected curve in $\PP^2_\Z$
of degree $3d$, hence of genus at least $1$,
so its Jacobian contradicts the main theorem of~\cite{fontaine1985}.
\end{example}

Despite these counterexamples,
P. Autissier has proved a positive result for a slightly different problem.
An {\em arithmetic variety} of dimension $m$
is an integral scheme $X$ of dimension $m$
that is projective and flat over $\Z$,
such that $X_\Q$ is regular (of dimension $m-1$).
If $\OO_K$ is the ring of integers of a finite extension $K$ of $\Q$,
then an {\em arithmetic variety over $\OO_K$} is an $\OO_K$-scheme $X$
such that $X$ is an arithmetic variety and whose generic fiber $X_K$
is geometrically irreducible over $K$.
The following is a part of Th\'eor\`eme~3.2.3 of~\cite{autissier2001}:

\begin{quote}
Let $X$ be an arithmetic variety over $\OO_K$ of dimension $m \ge 3$.
Then there exists a finite extension $L$ of $K$
and a closed subscheme $X'$ of $X_{\OO_L}$ such that
\begin{enumerate}
\item The subscheme $X'$ is an arithmetic variety over $\OO_L$
of dimension $m-1$.
\item Whenever the fiber $X_\pp$ of $X$ above $\pp \in \Spec \OO_K$
is smooth, the fiber $X'_{\pp'}$ of $X'$ above $\pp'$ is smooth
for all $\pp' \in \Spec \OO_L$ lying above $\pp$.
\end{enumerate}
\end{quote}
Actually Autissier proves more, that one can also control the height of $X'$.
(He uses the theory of heights developed by Bost, Gillet, and Soul\'e,
generalizing Arakelov's theory.)

The most significant difference between Autissier's result
and the phenomenon exhibited by Fakhruddin's examples 
is the finite extension of the base allowed in the former.

\section*{Acknowledgements}

I thank Ernie Croot for a conversation at an early stage of this research,
David Eisenbud for conversations about 
Section~\ref{largecharacteristicsection},
and Matt Baker and Najmuddin Fakhruddin for some other comments.
I thank Pascal Autissier and Ofer Gabber
for sharing their preprints with me,
and Jean-Pierre Jouanolou for sharing some unpublished notes about
resultants and discriminants.

\providecommand{\bysame}{\leavevmode\hbox to3em{\hrulefill}\thinspace}

\end{document}